\pdfoutput=1
\documentclass[preprint,3p,times]{elsarticle}

\usepackage{amsmath}
\usepackage[T1]{fontenc}
\usepackage{amssymb}
\usepackage{booktabs}
\usepackage{xcolor}
\usepackage{url} 

\journal{Journal of Computational and Applied Mathematics}

\renewcommand{\vec}[1]{\boldsymbol{#1}}

\usepackage[unicode=true, bookmarks=true, bookmarksnumbered=false, bookmarksopen=false, breaklinks=true, pdfborder={0 0 1}, colorlinks=true, linkcolor=black, urlcolor=theblue, citecolor=theblue, hyperfootnotes=false]{hyperref}
\hypersetup{pdftitle={New fully symmetric and rotationally symmetric cubature rules on the triangle using minimal orthonormal bases}, pdfauthor={Stefanos-Aldo Papanicolopulos}}

\begin{document}

\begin{frontmatter}

\title{New fully symmetric and rotationally symmetric cubature rules\\ on the triangle using minimal orthonormal bases}

\author[uoe]{Stefanos-Aldo~Papanicolopulos\corref{cor1}}
\ead{S.Papanicolopulos@ed.ac.uk}
\cortext[cor1]{Corresponding author.  Tel.: +44 (0)131 650 7214; Fax: +44 (0)131 650 6554.}
\address[uoe]{%
Institute for Infrastructure \& Environment, %
School of Engineering, The University of Edinburgh, %
Edinburgh, EH9 3JL, UK}

\begin{abstract}
Cubature rules on the triangle have been extensively studied, as they are of great practical interest in numerical analysis. In most cases, the process by which new rules are obtained does not preclude the existence of similar rules with better characteristics. There is therefore clear interest in searching for better cubature rules.  

Here we present a number of new cubature rules on the triangle, exhibiting full or rotational symmetry, that improve on those available in the literature either in terms of number of points or in terms of quality. These rules were obtained by determining and implementing minimal orthonormal polynomial bases that can express the symmetries of the cubature rules. As shown in specific benchmark examples, this results in significantly better performance of the employed algorithm.
\end{abstract}

\begin{keyword}
Cubature \sep triangle \sep fully symmetric rules \sep rotationally symmetric rules \sep symmetric polynomials
\MSC Primary 65D32 \sep Secondary 65D30
\end{keyword}

\end{frontmatter}

% ----------------------------------------------------------------------
\section{Introduction}

Cubature, that is the numerical computation of a multiple integral, is an important method of numerical analysis, as it is of great practical interest in different applications involving integration. An extensive literature therefore exists on this topic \citep[see e.g.][]{Stroud1971, Cools1997}, including also compilations of specific cubature rules \citep{Cools2003}.

The present paper considers cubature rules on the triangle. This is perhaps the most studied cubature domain, with a correspondingly large body of literature a selection of which is presented here. While rules of degree up to 20, thus covering most cases of practical interest, were progressively developed by 1985 \citep{Stroud1971, Cowper1973, LynessJespersen1975, Dunavant1985}, this is still an active field \citep{BerntsenEspelid1990, HeoXu1999, WandzuraXiao2003, Rathod2007865, ZhangCuiLiu2009, XiaoGimbutas2010, Williams2014, witherdenvincent2014, Papanicolopulos2015camwa, Witherden2015}. This happens for two distinct reasons, the first being that different applications require different properties of the cubature rules; the previously cited work for example focuses only on fully symmetric rules (which are also the easier to determine), while only a few works consider rotationally symmetric \citep{Gatermann1988,XiaoGimbutas2010} or asymmetric \citep{Taylor2007, Taylor2008} rules. The second reason explaining the interest in researching new cubature rules is that almost all rules in the literature have been determined numerically using an iterative procedure, therefore there is the possibility that a ``better'' rule (matching some given requirements) may exist, for example one having fewer points (see \citep{LynessCools1995} for a lower bound on the number of points for given degree). For fully symmetric rules, the fact that the ``best'' existing rules for degree up to 14 have indeed the minimal possible number of points was recently proved using solutions based on algebraic solving \citep{Papanicolopulos2015camwa}.

In this paper we focus on the iterative algorithm for obtaining fully symmetric cubature rules on the triangle initially proposed by \citet{ZhangCuiLiu2009} and recently refined by \citet{Witherden2015}. A main feature of~\citep{Witherden2015} (which had already been used in \citep{Taylor2007,XiaoGimbutas2010}) is the use of an orthonormal basis instead of the typical monomial basis usually employed. Further improving upon this point, we describe here a minimal orthonormal basis for fully symmetric rules and then extend this basis to also cover the case of rules with only rotational symmetry. This results in a number of new rules that improve upon those found in the literature, especially for the rotationally symmetric case.

The structure of the paper is as follows: after this introduction, Section~\ref{sec:background} summarises the required theoretical background. Section~\ref{sec:orthonormal} presents orthonormal bases for the fully symmetric case, including minimal bases, in terms of the typically used orthonormal polynomials while Section~\ref{sec:symcoordbasis} presents the minimal basis in terms of symmetric polynomials. In Section~\ref{sec:rotsym} we extend the minimal basis to obtain a minimal basis for cubature rules with rotational symmetry. A summary of the numerical results is presented and discussed in Section~\ref{sec:results}, while the conclusions of the paper are stated in Section~\ref{sec:conclusions}.

% ----------------------------------------------------------------------
\section{Theoretical background}
\label{sec:background}

A cubature rule approximates the integral of a function $f$ on a domain $\Omega$ (normalised by the domain's area $A$) as the weighted sum of the function's value evaluated at a set of $n_k$ points 
$\vec{x}_i$,
\begin{equation}
  \label{eq:gencub}
  \sum_i^{n_k} w_i f(\vec{x}_i) \approx \frac{1}{A} \int_{\Omega} f(\vec{x}) \mathrm{d} \vec{x}
\end{equation}
The cubature rule is of polynomial degree $\phi$ if equation~(\ref{eq:gencub}) is exact for all polynomials of degree up to $\phi$ but not exact for at least one polynomial of degree $\phi+1$. 

Since equation~(\ref{eq:gencub}) is linear in the function $f$, we only need to ensure that it is exact for a basis of the polynomials of degree $\phi$. The simplest such basis in two dimensions is  the set of monomials $x^i y^j$ in the Cartesian coordinates $x$ and $y$ with $i+j \leq \phi$, but for the triangle another simple basis is the set of monomials $L_1^i L_2^j L_3^{\phi-i-j}$ expressed in term of the areal (or barycentric) coordinates $L_1$, $L_2$ and $L_3$ (with all exponents being non-negative).     

In two dimensions each point contributes three unknowns (two coordinates and a weight), therefore setting in equation~(\ref{eq:gencub}) $f$ as each of the basis polynomials for degree $\phi$  results in a polynomial system of $(\phi+1)(\phi+2)/2$ equations in $3 n_k$ variables. The solution of this system yields the cubature point coordinates and weights defining the cubature rule.

In the general (asymmetric) case it can be quite difficult to solve the above-mentioned system even for moderate values of $\phi$, therefore some symmetry condition is imposed on the cubature points to reduce the number of unknowns. As mentioned in the introduction, these symmetries may also be a requirement of the application being considered; on the triangle, for example, full symmetry ensures that the computed approximate value of the integral is independent of the order in which the vertices are numbered.

For cubature rules on the triangle, the most commonly used symmetry is full symmetry, where if a point with areal coordinates $(L_1,L_2,L_3)$ appears in the rule, then all points resulting from permutation of the areal coordinates also appear. Depending on the number of distinct values of the areal coordinates we therefore obtain different symmetry orbits (for 1, 2 or 3 distinct values we get orbits of type 0, 1 or 2 which have 1, 3 or 6 points and contribute 1, 2 or 3 unknowns to the system of equations); see \citep{Witherden2015} for a more detailed explanation of symmetry orbits. Full symmetry allows for a significant reduction in the number of unknowns (roughly by a factor of 6 for larger values of $\phi$) and, through appropriate considerations, for a corresponding decrease in the number of equations \citep{LynessJespersen1975}. 

The disadvantage of full symmetry is that in most cases it does not lead to the cubature rule with the minimal number of points for a given degree and quality (in the sense of ``quality'' defined in Section~\ref{sec:results}). It is possible to get rules with fewer points, while still reducing the number of equations and unknowns, by requiring only rotational symmetry. In this case instead of considering all the permutations of the areal coordinates we only consider the even permutations, so that $(L_1,L_2,L_3)$ is permuted into $(L_2,L_3,L_1)$ and $(L_3,L_1,L_2)$. This results in two types of orbits: type-0 with only one point (the centroid) and type-1 with three points, therefore the number of unknowns is approximately twice that of the fully symmetric case.

% ----------------------------------------------------------------------
\section{Orthonormal bases on the triangle}
\label{sec:orthonormal}

\subsection{A full orthonormal basis}

While the monomials (in either the Cartesian or the areal coordinates) described in Section~\ref{sec:background} are the simplest basis polynomials, they lead at higher degrees to polynomial systems which are poorly conditioned, therefore the use of an orthonormal basis has been proposed \citep{Taylor2007,XiaoGimbutas2010,Witherden2015}.

A standard set of orthonormal basis polynomials on the triangle has been proposed in the literature \citep{proriol1957,koornwinder1975,Dubiner1991}, which we can write in the form
\begin{equation}
  \psi_{ij}(\vec{x}) = \hat{P}_i \big( d/s \big) \hat{P}_j^{(2i+1,0)} \big( 1 - 2 s \big) s^i
\end{equation}
where $\hat{P}^{(\alpha,\beta)}_n = \sqrt{2n+\alpha+1} P^{(\alpha,\beta)}_n $ are scaled Jacobi polynomials and the values $d$ and $s$ depend on the coordinates.

Specific expressions for the $\psi_{ij}(\vec{x})$ (and therefore for $d$ and $s$ in Cartesian coordinates) are given in the literature by specifying a reference triangle. Using areal coordinates, however, $d$ and $s$ are simply the difference and sum respectively of two of the areal coordinates, without reference to a specific triangle. Choosing for example $L_2$ and $L_1$ we have
\begin{equation}
  s = L_2 + L_1, \quad d = L_2 - L_1
\end{equation}

An interesting property of the basis polynomials expressed in terms of $d$ and $s$ is that the $\psi_{ij}$ are the Gram-Schmidt orthonormalisation of the monomials $d^i (-s)^j$ taken in increasing graded lexicographic order.

\medskip

The basis $\mathcal{\tilde{P}}^\phi_{f}$ of all polynomials of degree up to $\phi$ (the ``full'' basis) then contains all basis polynomials with $i + j \leq \phi$, that is
\begin{equation}
  \label{eq:basis_full}
  \mathcal{\tilde{P}}^\phi_{f} = \bigl\{\psi_{ij}(\vec{x}) \mid 0 \leq
  i \leq \phi,\; 0 \leq j \leq \phi - i \bigr\}
\end{equation}
Indicating by $n(\phi)$ the cardinality of the degree-$\phi$ basis and by $m(\omega)$ the number of basis polynomials of degree $\omega$, for the full basis we easily see that
\begin{equation}
  \label{eq:card_full}
  n_f(\phi) = \frac{(\phi+1)(\phi+2)}{2} , \quad
  m_f(\omega) = \omega+1
\end{equation}

\subsection{Objective orthonormal bases for fully symmetric rules}
\label{sec:objsymbasis}

While a full basis is needed to represent all polynomials of degree $\phi$, a reduced basis can be used when considering fully symmetric cubature rules, as this restricts the form of the system of polynomial equations to be solved. \citet{Witherden2015} propose an ``objective'' basis, that is a subset of the full basis that can still represent the polynomial system for fully symmetric rules,
\begin{equation}
\label{eq:basistriWV}
\mathcal{\tilde{P}}^\phi_{w} = \bigl\{\psi_{ij}(\vec{x}) \mid
         0 \leq i \leq \phi,\; i \leq j \leq \phi - i \bigr\}
\end{equation}
In equation~(\ref{eq:basistriWV}), due to the limits on $j$, the actual limits on $i$ are $0 \leq i \leq \lfloor \phi/2 \rfloor$, and it is therefore easy to show that
\begin{equation}
  \label{eq:card_wv}
  n_w(\phi) = \left\lfloor \frac{(\phi+2)^2}{4} \right\rfloor \sim \frac{1}{2} n_f(\phi) , \quad
  m_w(\omega) = 1 + \lfloor \omega/2 \rfloor	
\end{equation}
\ref{sec:cardinality} presents the derivation of equations~(\ref{eq:card_wv}), with results for the bases given below being obtained in a similar way.

As already noted in \citep{Witherden2015}, this objective basis is not optimal as its modes are not completely independent. Indeed, while the basis $\mathcal{\tilde{P}}^\phi_{w}$ is an objective basis, it is interesting to note that there is no obvious reason why the specific $\psi_{ij}$ polynomials were omitted. It is actually possible to have other objective bases with the same cardinality that use another subset of the full basis, such as
\begin{equation}
\label{eq:basistriWVstar}
\mathcal{\tilde{P}}^\phi_{w_2} = \bigl\{\psi_{ij}(\vec{x}) \mid
         0 \leq j \leq \lfloor \phi/2 \rfloor,\; j \leq i \leq \phi - j \bigr\}
\end{equation}
which is actually $\mathcal{\tilde{P}}^\phi_{w_2}$ with the indices $i$ and $j$ swapped.

To further reduce the cardinality of the basis, we first note that for a symmetric orbit we will be adding the polynomials $\psi_{ij}(d,s)$ and $\psi_{ij}(-d,s)$, which correspond to points with areal coordinates $(L_1,L_2,L_3)$ and $(L_2,L_1,L_3)$. If $i$ is odd, however, $\psi_{ij}(d,s)$ is also odd with respect to d and therefore all $\psi_{ij}$ with odd $i$ can be removed from the objective basis $\mathcal{\tilde{P}}^\phi_{w}$ to obtain the ``even'' basis%
\footnote{Obviously different even bases can be obtained, e.g. starting from the basis  $\mathcal{\tilde{P}}^\phi_{w_2}$.}
\begin{equation}
\label{eq:basistrieven}
\mathcal{\tilde{P}}^\phi_{e} = \bigl\{\psi_{ij}(\vec{x}) \mid 0 \leq
  i \leq \lfloor\phi/2\rfloor,\; i \leq j \leq \phi - i, i\;\mathrm{even} \bigr\}
\end{equation}
for which we obtain
\begin{equation}
  n_e(\phi) = \left\lfloor \frac{(\phi+3)^2}{8} \right\rfloor  \sim \frac{1}{4} n_f(\phi) ,\quad
  m_e(\omega) = 1 + \lfloor \omega/4 \rfloor
\end{equation}

Expressing the basis polynomials in terms of $d$ and $s$, and then $d$ and $s$ in terms of the areal coordinates, has therefore the advantage of making obvious the symmetry and antisymmetry of the basis polynomials with respect to exchange of two vertices.

\medskip

As will also be discussed in Section~\ref{sec:symcoordbasis}, a symmetric basis with even lower cardinality is possible. Indeed for the minimal basis we get~\citep{LynessCools1995}
\begin{equation}
  \label{eq:nm_min}
  n_m(\phi) = \left\lfloor \frac{(\phi+3)^2}{12} + \frac{1}{4} \right\rfloor
      \sim \frac{1}{6} n_f(\phi) ,\quad
  m_m(\omega) = 1 + \lfloor \omega/6 \rfloor - \kappa_6(\omega)
\end{equation}
where
\begin{equation}
  \label{eq:kappa_a}
  \kappa_a(\omega) = \begin{cases} 1 & \textrm{if $\omega \bmod a = 1$} \\ 0 & \textrm{otherwise} \end{cases}
\end{equation}

To construct a minimal basis of degree $\phi$, it suffices to choose a subset of the even base of the same degree $\phi$ so that the number of polynomials $\psi_{ij}$ with $i+j=\omega$ is given by $m_m(\omega)$ as defined in equation~(\ref{eq:nm_min}), that is
\begin{equation}
\label{eq:basistrimingen}
\mathcal{\tilde{P}}^\phi_{\bar{m}} = \bigl\{\psi_{ij}(\vec{x}) \in \mathcal{\tilde{P}}^\phi_{e} \mid
  \# \{ \psi_{ij} \mid i+j=\omega \leq \phi \} = m_m(\omega) \bigr\}
\end{equation}
While we do not provide here a proof that $\mathcal{\tilde{P}}^\phi_{\bar{m}}$ is indeed an objective basis, it is relatively easy to check this for given values of $\phi$ using a computer algebra system.

It is easy to create two such minimal objective bases as
\begin{gather}
\label{eq:basistrimin_alt}
  \mathcal{\tilde{P}}^\phi_{m} = \bigl\{\psi_{2i,\omega-2i}(\vec{x})
    \mid 0 \leq i \leq m_m(\omega)-1,\; 0 \leq \omega \leq \phi \bigr\}
  \\
\label{eq:basistrimin_alt2}
  \mathcal{\tilde{P}}^\phi_{m_2} = \bigl\{\psi_{2i,\omega-2i}(\vec{x})
    \mid \lfloor \omega/2 \rfloor - \big( m_m(\omega)-1 \big) \leq i \leq \lfloor \omega/2 \rfloor,\;
         0 \leq \omega \leq \phi \bigr\}
\end{gather}
These two bases are obtained by considering each polynomial degree $0 \leq \omega \leq \phi$ and taking $m_m(\omega)$ consecutive even values for the first index of the basis polynomials, with the second index defined by the requirement that the sum of the two indices is equal to the degree $\omega$; in the first case we get the lowest possible values for the first index, while in the second case we get the highest possible values.

The bases~(\ref{eq:basistrimin_alt}) and~(\ref{eq:basistrimin_alt2}) can alternatively be written, after some calculations, as
\begin{gather}
\label{eq:basistrimin}
  \mathcal{\tilde{P}}^\phi_{m} = \bigl\{\psi_{ij}(\vec{x}) \mid 0 \leq
  i \leq \lfloor\phi/3\rfloor,\; 2i \leq j \leq \phi - i, i\;\mathrm{even}, j \neq 2i+1 \bigr\}
   \\
\label{eq:basistrimin2}
  \mathcal{\tilde{P}}^\phi_{m_2} = \bigl\{\psi_{ij}(\vec{x}) \mid
     0 \leq i \leq \phi,\; 0 \leq j \leq \min( \phi - i, i/2 ), i\;\mathrm{even} \bigr\}
\end{gather}
This form, though less intuitive, was found to be slightly simpler to implement in a computer code.

Other minimal objective bases can also be derived. For example, minimising the maximum of $i$ and $j$ for a given $\phi$ (trying to reduce number of computations and round-off error), results in the following basis
\begin{equation}
  \label{eq:basistrimin3}
  \mathcal{\tilde{P}}^\phi_{m_3} = \bigl\{\psi_{ij}(\vec{x}) \mid
    0 \leq i \leq 2 \lfloor\phi/3\rfloor + 2 \kappa_6(\phi-1),\;
    2 \lfloor i/4 \rfloor \leq j \leq \min\{\phi - i,2i\},\; i\;\mathrm{even} \bigr\}
\end{equation}

\begin{figure}
  \centering
  \includegraphics[width=\textwidth]{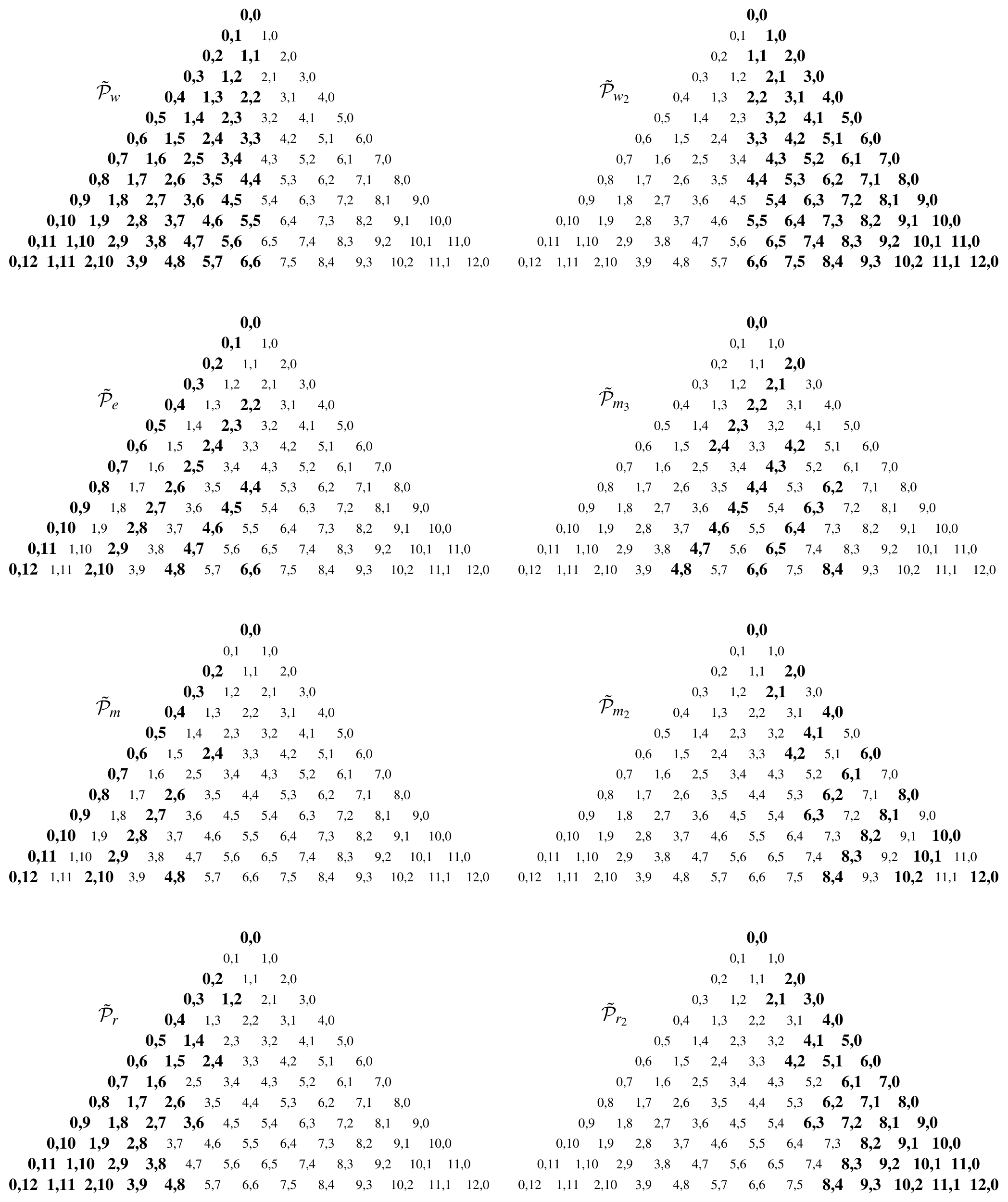}
  \caption{Basis polynomials used in the different objective bases presented in this paper, for polynomial degree $\phi=12$. Only the index pairs $i,j$ are shown instead of $\psi_{ij}$; the bold font indicates basis polynomials that are used in the objective basis, while smaller font size indicates polynomials in the full basis that are not used in the objective basis.}
  \label{fig:bases}
\end{figure}

Figure~\ref{fig:bases} provides a graphical representation of the objective bases presented in this paper (including those for rotational symmetry, presented in section~\ref{sec:rotsym}), for polynomial degree $\phi=12$. This should provide a more intuitive understanding of how each objective basis is obtained. As the rows of each pyramid correspond to basis polynomials of equal degree, the graphical representations for lower values of $\phi$ are obtained by removing rows from the bottom of each pyramid.

% ----------------------------------------------------------------------
\section{An orthonormal basis for symmetric polynomials}
\label{sec:symcoordbasis}

We present here the derivation of an orthonormal basis to be used in computing fully symmetric cubature rules on the triangle, which makes full use of the imposed symmetry by using symmetric polynomials \citep{Macdonald1998}. More details on this approach can be found in \citep{Papanicolopulos2015camwa}.

A \emph{symmetric polynomial} is a multivariate polynomial in $n$ variables, say $x_1,x_2,\ldots,x_n$, which is invariant under any permutation of its variables. We define the elementary symmetric polynomials $\tilde{x}_k$ as the sums of all products of $k$ distinct variables $x_i$, with negative sign when $k$ is odd, that is
\begin{equation}
  \label{eq:elemsym}
  \tilde{x}_k = (-1)^k \sum_{i_1<i_2<\ldots<i_k} x_{i_1} x_{i_2} \cdots x_{i_k}
\end{equation}
with $\tilde{x}_0 = 1$.
The fundamental theorem of symmetric polynomials states that any symmetric polynomial in the variables $x_i$ can be expressed as a polynomial in the elementary symmetric polynomials $\tilde{x}_k$.

We consider, without loss of generality, a type-2 orbit in a fully symmetric cubature rule on the triangle. This orbit consists of a point with areal coordinates $(L_1,L_2,L_3)$ and the five other points resulting from permutation of these coordinates (which, for the type-2 orbit, are all distinct). Using equation~(\ref{eq:gencub}) for a polynomial $\hat{f}(L_1,L_2,L_3)$ in the areal coordinates, yields only sums of the form%
\footnote{These sums appear multiplied by the weight corresponding to the orbit being considered.}
\begin{equation}
  \label{eq:Ts}
  T_s =
  \hat{f}(L_1,L_2,L_3) + \hat{f}(L_3,L_1,L_2) + \hat{f}(L_2,L_3,L_1) + 
  \hat{f}(L_1,L_3,L_2) + \hat{f}(L_2,L_1,L_3) + \hat{f}(L_3,L_2,L_1) 
\end{equation}

Therefore, for fully symmetric rules, the left hand side of~(\ref{eq:gencub}) only contains symmetric polynomials in the areal coordinates. According to the fundamental theorem of symmetric polynomials, these can therefore be expressed as polynomials in the elementary symmetric polynomials $\tilde{L}_1 = -(L_1+L_2+L_3)$, $\tilde{L}_2 = L_1 L_2 + L_2 L_3 + L_3 L_1$ and $\tilde{L}_3 = - L_1 L_2 L_3$ (noting, however, that in this case $\tilde{L}_1=-1$).

It is therefore easily seen that instead of considering all polynomials of degree $\phi$, or at least all polynomials in a basis of degree $\phi$, we only need to consider a symmetric basis consisting of the largest possible number of linearly independent polynomials in $\tilde{L}_2$ and $\tilde{L}_3$ of weighted total degree less or equal to $\phi$ (with a weight 2 for $\tilde{L}_2$ and a weight 3 for $\tilde{L}_3$, as they respectively involve double and triple products of $L_1$, $L_2$, and $L_3$).

The simplest such symmetric basis consists of the monomials $\tilde{L}_2^i \tilde{L}_3^j$ with $2i + 3j \leq \phi$, that is
\begin{equation}
  \label{eq:symbasismonom}
  \mathcal{Q}^\phi_{s} = \bigl\{ \tilde{L}_2^i \tilde{L}_3^j \mid 2i + 3j \leq \phi \bigr\}
\end{equation}

For the basis $\mathcal{Q}^\phi_{s}$ (and indeed for any symmetric basis) we easily obtain
\begin{equation}
  \label{eq:nm_sym}
  n_s(\phi) = \left\lfloor \frac{(\phi+3)^2}{12} + \frac{1}{4} \right\rfloor
      \sim \frac{1}{6} n_f(\phi) ,\quad
  m_s(\omega) = 1 + \lfloor \omega/6 \rfloor - \kappa_6(\omega)
\end{equation}
with $\kappa_a(\omega)$ already defined in equation~(\ref{eq:kappa_a}), since we already used this result in Section~\ref{sec:objsymbasis} for the minimal objective basis.

The monomial symmetric basis given in equation~(\ref{eq:symbasismonom}) is obviously not orthogonal. To obtain an orthonormal symmetric basis $\mathcal{\tilde{Q}}^\phi_{s}$ we can orthonormalise the monomials in the basis $\mathcal{Q}^\phi_{s}$. While the orthonormalisation can be done numerically, to minimise numerical errors we choose here to perform it analytically with a computer algebra system. This also allows for an efficient implementation of a multivariate Horner scheme \citep{Pena99}. Note that monomials must be considered in weighted lexicographic order to obtain orthonormal bases which include the bases of lower degree. \Citet{Chabysheva2013} have recently discussed an orthonormalisation of this type, but their use of Cartesian coordinates leads to polynomials with a significantly larger number of terms, and of higher degree.

It is important to note that the minimal objective basis and the orthonormal symmetric basis are not bases of the same polynomials. Indeed, the minimal objective basis is not a proper basis of the symmetric polynomials (and actually does not consist of symmetric polynomials); we need to sum the values of the basis polynomials $\psi_{ij}$ on all points of the orbit to obtain a basis for the symmetric polynomials (which is then no-longer orthogonal). The orthonormal symmetric basis, on the other hand, is a proper orthonormal basis of the symmetric polynomials.

It is not clear whether the fact that the symmetric basis is really orthonormal would by itself provide better efficiency or accuracy in obtaining results; the obvious advantage of the symmetric orthonormal basis is that it requires only a single evaluation of the basis polynomials instead of the six evaluations (for type-2 orbits) required by the minimal objective basis. On the other hand, the advantage of the minimal objective basis is that it is expressed in analytical form (in terms of the Jacobi polynomials), making it easier to implement in a computer code. Additionally, the product form of the polynomials in the objective basis allow for their more efficient evaluation.

Using either type of basis will result in a polynomial system with solutions that correspond to the same set of cubature rules, but only if complex solutions are also taken into account. It is of theoretical interest that there can be real solutions of the polynomial system expressed in terms of the symmetric polynomials that correspond to cubature rules with real weights but complex point coordinates. Using an objective basis, on the other hand, it is obvious that any real solution corresponds to a rule with real weights and coordinates.

% ----------------------------------------------------------------------
\section{Rotational symmetry}
\label{sec:rotsym}

We consider now a rule with rotational symmetry. The system of polynomial equations will now contain, instead of the terms $T_s$ in equation~(\ref{eq:Tr}), polynomials in the areal coordinates of the form
\begin{equation}
  \label{eq:Tr}
  T_r =
  \hat{f}(L_1,L_2,L_3) + \hat{f}(L_3,L_1,L_2) + \hat{f}(L_2,L_3,L_1) 
\end{equation}
This can be written as
\begin{equation}
  T_r = \frac{ T_s + T_a }{2}
\end{equation}
where $T_s$ is the symmetric polynomial given in equation~(\ref{eq:Ts}) and $T_a$ is the antisymmetric polynomial
\begin{equation}
  \label{eq:Ta}
  T_a =
  \hat{f}(L_1,L_2,L_3) + \hat{f}(L_3,L_1,L_2) + \hat{f}(L_2,L_3,L_1) - 
  \hat{f}(L_1,L_3,L_2) - \hat{f}(L_2,L_1,L_3) - \hat{f}(L_3,L_2,L_1) 
\end{equation}

As already mentioned, $T_s$ can be expressed as a polynomial in the symmetric polynomials $\tilde{L}_2$ and $\tilde{L}_3$. The antisymmetric polynomial, on the other hand, can be expressed as the product of a symmetric polynomial (in $\tilde{L}_2$ and $\tilde{L}_3$) with the alternating polynomial $\tilde{L}_A$,
\begin{equation}
  \tilde{L}_A = (L_1-L_2)(L_1-L_3)(L_2-L_3)
\end{equation}
Considering that $\tilde{L}_A$ is of degree 3 in the areal coordinates, we see that a rotationally symmetric basis using monomials is given by
\begin{equation}
  \label{eq:rotbasismonom}
  \mathcal{Q}^\phi_{r} = \bigl\{ \tilde{L}_2^i \tilde{L}_3^j \tilde{L}_A^k  \mid 2i + 3j + 3k \leq \phi,\; k \in \{0,1\} \bigr\}
\end{equation}
from which we can obtain~\citep{LynessCools1995}
\begin{equation}
  \label{eq:nm_rot}
  n_r(\phi) = 1 + \left\lfloor \frac{(\phi+3)\phi}{6} \right\rfloor
      \sim \frac{1}{3} n_f(\phi) ,\quad
  m_r(\omega) = 1 + \lfloor \omega/3 \rfloor - \kappa_3(\omega)
\end{equation}

As in the fully symmetric case, the monomials in the basis $\mathcal{Q}^\phi_{r}$ can be orthonormalised to obtain an orthonormal rotationally symmetric basis $\mathcal{\tilde{Q}}^\phi_{r}$. 

It is also possible to obtain minimal objective bases for rotationally symmetric rules in terms of the basis polynomials $\psi_{ij}$. After some calculations it can be seen that these will consist of a minimal objective basis for fully symmetric rules plus a set of basis polynomials $\psi_{ij}$ with $i$ odd. The bases in equations~(\ref{eq:basistrimin}) and~(\ref{eq:basistrimin2}) yield the following minimal objective bases for rotationally symmetric rules
\begin{gather}
\label{eq:basistriminrot}
  \mathcal{\tilde{P}}^\phi_{r} = \bigl\{\psi_{ij}(\vec{x}) \mid 0 \leq
  i \leq \lfloor\phi/3\rfloor,\; 2i \leq j \leq \phi - i, j \neq 2i+1 \bigr\}
   \\
\label{eq:basistriminrot2}
  \mathcal{\tilde{P}}^\phi_{r_2} = \bigl\{\psi_{ij}(\vec{x}) \mid
     0 \leq i \leq \phi,\; 0 \leq j \leq \min( \phi - i, \lfloor i/2 \rfloor - \kappa_2(i) ) \bigr\}
\end{gather}
As can be seen in Figure~\ref{fig:bases}, $\mathcal{\tilde{P}}^\phi_{r_2}$ is just $\mathcal{\tilde{P}}^\phi_{r}$ with the indices $i$ and $j$ swapped.

% ----------------------------------------------------------------------
\section{Results}
\label{sec:results}

\subsection{Performance measurements}

\Citet{Witherden2015} have develop the C++ code \texttt{polyquad} to compute fully symmetric cubature rules (on the triangle and on other domains) using objective orthonormal bases. The objective bases for fully or rotationally symmetric rules proposed in Section~\ref{sec:orthonormal} can be easily implemented with minor modifications to the existing \texttt{polyquad} code.
The orthonormal basis for symmetric and rotationally symmetric polynomials presented in Section~\ref{sec:symcoordbasis} could also be implemented in \texttt{polyquad}, requiring however more extensive changes to the code. It was therefore found simpler to implement the algorithm in a new Fortran 95 code named \texttt{pq}. The two implementations are not directly comparable, and their relative performance depends among others on the minimisation solver used and its parameters. Comparing the two codes does however provide a first insight on the ability of one method to outperform the other.

Table~\ref{tbl:performance} shows the performance of \texttt{polyquad} and \texttt{pq} for the case of fully symmetric rules of degree $\phi=15$ with $49$ points, considering four different combination of orbits (using the notation $[i,j,k]$ to indicate a rule with $i$ type~0 orbits, $j$ type~1 orbits and $k$ type~2 orbits), of which only the combinations $[1,4,6]$ and $[1,6,5]$ actually yield a cubature rule. The performance is expressed as the number of trial rules evaluated per second, and represent the average of 20 different runs with at least 100 rules evaluated per run.

Similarly, table~\ref{tbl:performanceR} shows the performance of \texttt{polyquad} and \texttt{pq} for the case of rotationally symmetric rules. As in this case there is a single combination of orbits for a given total number of points, rules of different degrees, from 13 to 16, were evaluated so as to always consider a combination of degree and number of points that actually yields cubature rules.

\begin{table}
\caption{Performance (rules/sec) for degree-15 rules with 49 points, using $\texttt{polyquad}$ (with different objective bases) and using \texttt{pq}.}
\label{tbl:performance}
\center
\begin{tabular}{ccccccccc}
\toprule
rule & $\mathcal{\tilde{P}}^{15}_{f}$ & $\mathcal{\tilde{P}}^{15}_{w}$  & $\mathcal{\tilde{P}}^{15}_{w_2}$ & $\mathcal{\tilde{P}}^{15}_{e}$ & $\mathcal{\tilde{P}}^{15}_{m}$ & $\mathcal{\tilde{P}}^{15}_{m_2}$ & $\mathcal{\tilde{P}}^{15}_{m_3}$ & $\mathcal{\tilde{Q}}^{15}_{s}$ \\
\midrule
{}[1,2,7] & 1.06 & 1.80 & 1.64 & 2.85 & 3.70 & 3.72 & 2.89 & 6.26 \\
{}[1,4,6] & 1.00 & 1.77 & 1.55 & 2.59 & 3.18 & 3.52 & 2.61 & 4.63 \\
{}[1,6,5] & 1.09 & 1.88 & 1.74 & 2.89 & 3.58 & 3.63 & 2.88 & 4.68 \\
{}[1,8,4] & 1.71 & 2.79 & 2.41 & 4.20 & 4.79 & 5.00 & 3.67 & 4.35 \\
\bottomrule
\end{tabular}
\end{table}

\begin{table}
\caption{Performance (rules/sec) for four different types of rules with rotational symmetry, using $\texttt{polyquad}$ with two different objective bases and using \texttt{pq}.}
\label{tbl:performanceR}
\center
\begin{tabular}{ccccc}
\toprule
degree & points & $\mathcal{\tilde{P}}^{\phi}_{r}$ & $\mathcal{\tilde{P}}^{\phi}_{r_2}$  & $\mathcal{\tilde{Q}}^{\phi}_{r}$ \\
\midrule
13 & 36 & 2.66 & 2.84 & 2.80 \\
14 & 42 & 1.23 & 1.25 & 1.61 \\
15 & 46 & 0.92 & 0.97 & 1.30 \\
16 & 52 & 0.56 & 0.59 & 0.66 \\
\bottomrule
\end{tabular}
\end{table}

While the exact values depend on the compiler and hardware used, the results in tables~\ref{tbl:performance} and~\ref{tbl:performanceR} show the relative performance of different bases, with the best results in \texttt{polyquad} obtained with the minimal objective bases $\mathcal{\tilde{P}}_{m_2}$ and $\mathcal{\tilde{P}}_{r_2}$, and with the use of orthonormal bases in \texttt{pq} outperforming the use of objective bases in \texttt{polyquad}.

The performance obtained using the orthonormal basis for symmetric polynomials critically depends on the efficiency with which the basis polynomials can be evaluated. It is actually expected that appropriate optimisation of the computation of the objective basis could lead to faster evaluation than in the case of the symmetric basis.

\subsection{New cubature rules}

In presenting specific rules we are interested in the ``quality'' of the rule, which is expressed using two letters. The first letter is `P' if all weights are positive (otherwise it's `N') and the second is `I' if all points of the rule lie within the triangle (otherwise it's `O'). We therefore obtain PI, NI, PO and NO rules (in decreasing order of quality). 

Table~\ref{tbl:fulsym} presents the newly obtained fully symmetric rules that improve on the ones in the literature either on quality or on number of points. As mentioned in the introduction, it is already known~\cite{Papanicolopulos2015camwa} that no improved rules could be found for $\phi \leq 14$. For PI rules, no improved results were obtained for $\phi \leq 23$. 

\begin{table}
\caption{Number of points and quality for new fully symmetric cubature rules}
\label{tbl:fulsym}
\center
\begin{tabular}{crrrrrrrrrrr}
\toprule
 & \multicolumn{11}{c}{degree} \\
\cmidrule{2-12}
quality & 15 & 16 & 17 & 18 & 19 & 20 & 21 & 22 &  23 &  24 &  25 \\
\midrule  
PI      &    &    &    &    &    &    &    &    & 102 &     &     \\
NI      & 48 &    &    &    &    &    &    &    &     &     &     \\ 
PO      &    &    & 58 &    &    &    &    & 94 &     &     & 118 \\
\bottomrule
\end{tabular}
\end{table}

The fully symmetric case has been extensively studied in the literature, especially for degrees up to 20, therefore only a few new results were found. The implementation of the rotationally symmetric basis, on the other hand, yielded a larger number of new rules that improve in some way on the results previously available (either in number of points or in quality for a given number of points). These new rules are summarised in Table~\ref{tbl:rotsym}, starting from degree 12 as for lower degrees no improved rules were obtained. Many of the obtained rules are of PI quality. However, when NI (or PO) rules with fewer points were encountered these are also mentioned.
\begin{table}
\caption{Number of points and quality for new rotationally symmetric cubature rules}
\label{tbl:rotsym}
\center
\begin{tabular}{crrrrrrrrrrrrrr}
\toprule
 & \multicolumn{14}{c}{degree} \\
\cmidrule{2-15}
quality & 12 & 13 & 14 & 15 & 16 & 17 & 18 & 19 & 20 & 21 & 22 &  23 &  24 &  25 \\
\midrule
PI      &    &    &    &    &    & 57 &    & 70 &    & 85 &    & 100 & 109 & 117 \\
NI      &    &    &    &    &    &    & 64 &    &    &    &    &     &     &     \\ 
PO      & 31 &    & 40 &    & 51 &    &    &    &    &    &    &     &     &     \\
\bottomrule
\end{tabular}
\end{table}

The coordinates and weights for the rules summarised in Tables~\ref{tbl:fulsym} and~\ref{tbl:rotsym}, computed to double precision, are provided as ancillary files at \url{http://arxiv.org/abs/1411.5631v2}. While in most cases many rules were computed for given degree, number of points and quality, only one rule of each type is presented. This rule was selected to minimise the ratio of maximum to minimum weight, avoiding however (for PI rules) rules with points almost on the boundary.

\medskip

Rules of increasing degree take longer to be computed, and are of decreasing interest in practical applications. For this reason, only rules of degree up to 25 have been considered here. There is however no indication that rules of higher degree cannot be obtained using the same method, given enough computation time. It is on the other hand also possible that improved rules may be obtained even for the degrees considered here.

% ----------------------------------------------------------------------
\section{Conclusions}
\label{sec:conclusions}

We have presented in this paper minimal orthonormal polynomial bases on the triangle for computing fully symmetric and rotationally symmetric cubature rules. These bases can be either ``objective'' bases, that is subsets of the complete polynomial basis that yield the required symmetry, or true fully/rotationally symmetric bases in terms of the symmetric elementary polynomials (and the alternating polynomial for rotational symmetry).

As these bases are minimal, they allow for more efficient  computation of cubature rules. We therefore present a number of new rules that improve, in some aspects, on the rules available in the literature. Especially for the rotationally symmetric rules, a large number of new rules is obtained, most of which are of PI quality.

Further optimisation of the implementation of the algorithm could be possible, for example by implementing a more efficient computation of the basis polynomials or by employing a different optimisation solver to solve the polynomial equations. This is currently a work in progress, as it would allow more efficient computation of rules of higher degree, should they be needed, and especially more efficient computation of cubature rules on the tetrahedron.

\section*{Acknowledgements}

This research effort is funded from the People Programme (Marie Curie Actions) of the European Union's Seventh Framework Programme (FP7/2007-2013) under REA grant agreement n\textsuperscript{o} 618096.

\appendix

\section{Cardinality of the objective bases}
\label{sec:cardinality}

To compute the cardinality of the objective bases presented in this paper, we make use of the well known formula
\begin{equation}
  \sum_{i=0}^{\nu} i = \frac{\nu(\nu+1)}{2}
\end{equation}

Considering first the full basis, we recall here equation~(\ref{eq:basis_full})
\begin{equation*}
  \mathcal{\tilde{P}}^\phi_{f} = \bigl\{\psi_{ij}(\vec{x}) \mid 0 \leq
  i \leq \phi,\; 0 \leq j \leq \phi - i \bigr\}
\end{equation*}
For each value of $i$ we have $\phi-i+1$ values of $j$. The cardinality of the basis is therefore
\begin{equation}
  n_f(\phi) = \sum_{i=0}^{\phi} (\phi-i+1)
            = \sum_{i=0}^{\phi} (\phi+1) - \sum_{i=0}^{\phi} i
            = (\phi+1)(\phi+1) - \frac{\phi(\phi+1)}{2}
            = \frac{(\phi+2)(\phi+1)}{2}
\end{equation}
which is the well known result given in equation~(\ref{eq:card_full}). The number of basis polynomials of degree $\omega$ is then directly calculated as
\begin{equation}
  m_f(\omega) = n_f(\omega) - n_f(\omega-1) = \omega + 1
\end{equation}

The same results can be obtained by calculating first $m_f(\omega)$. Considering that the degree $\omega$ of the basis polynomial $\psi_{ij}$ is simply $\omega=i+j$, we replace $j=\omega-i$ in the inequalities
\begin{equation}
  0 \leq i \leq \phi,\; 0 \leq j \leq \phi - i
\end{equation}
to obtain, after some very simple manipulations,
\begin{equation}
  0 \leq \omega \leq \phi ,\; 0 \leq i \leq \omega
\end{equation}
From the second set of inequalities we directly obtain $m_f(\omega) = \omega + 1$, so the cardinality is easily computed as
\begin{equation}
  n_f(\phi) = \sum_{\omega=0}^{\phi} (\omega+1) = \frac{(\phi+2)(\phi+1)}{2}
\end{equation}

Either of the procedures described above for the full basis can be employed to obtain the cardinality of the other objective bases. The main aspect to consider is that the bounds for $i$ and $j$ (or $\omega$) must be strict. This introduces the need to use the floor function, which in turn makes the computations slightly more complicated.

Consider for example the objective basis introcuded in \citep{Witherden2015}, given here as $\mathcal{\tilde{P}}^\phi_{w}$ in equation~(\ref{eq:basistriWV}). The bounds for $i$ and $j$ are given by
\begin{equation}
  0 \leq i \leq \phi,\; i \leq j \leq \phi - i
\end{equation}
From the first and last term in the second set of inequalities we get $2i \leq \phi$. As $\phi$ can be odd, the strict bound for $i$ is $i \leq \lfloor \phi/2 \rfloor$ so that the strict bounds are
\begin{equation}
  \label{eq:ineq_ij_wv}
  0 \leq i \leq \lfloor\phi/2 \rfloor,\; i \leq j \leq \phi - i
\end{equation}
For each value of $i$ we have $\phi-2i+1$ values of $j$ therefore the cardinality of the basis is
\begin{equation}
  n_w(\phi) = \sum_{i=0}^{\lfloor\phi/2 \rfloor} (\phi-2i+1)
            = \sum_{i=0}^{\lfloor\phi/2 \rfloor} (\phi+1) - 2 \sum_{i=0}^{\lfloor\phi/2 \rfloor} i
            = (\phi+1)(\lfloor\phi/2 \rfloor + 1) - 2 \frac{\lfloor\phi/2 \rfloor(\lfloor\phi/2 \rfloor+1)}{2}
            = (\phi+1 - \lfloor\phi/2 \rfloor)(\lfloor\phi/2 \rfloor+1)
\end{equation}
Considering $\phi$ as either even or odd we can write $\phi=2k+l$ with $l \in \{0,1\}$ so that $\lfloor\phi/2 \rfloor = (\phi-l)/2$ and
\begin{equation}
  n_w(\phi) = \big(\phi+1 - (\phi-l)/2\big)\big((\phi-l)/2+1\big)
            = \frac{(\phi+2)^2}{4} - \frac{l^2}{4}
            = \left\lfloor \frac{(\phi+2)^2}{4} \right\rfloor
\end{equation}
where the last step is computed considering that $l^2/4 < 1$ and that $n_w$ is integer.
We can then calculate $m_w$, setting now $\omega=2k+l$ with $l \in \{0,1\}$, as
\begin{align}
  m_w(\omega) &= n_w(\omega) - n_w(\omega-1)
               = \left\lfloor \frac{(\omega+2)^2}{4} \right\rfloor - \left\lfloor \frac{(\omega-1+2)^2}{4} \right\rfloor
               = \left\lfloor \frac{(2k+l+2)^2}{4} \right\rfloor - \left\lfloor \frac{(2k+l-1+2)^2}{4} \right\rfloor \notag\\
              &= \left\lfloor (k+1)^2 + (k+1)l + l^2/4 \right\rfloor - \left\lfloor (k+1)^2 + (k+1)(l-1) + (l-1)^2/4 \right\rfloor\notag\\
              &= (k+1)^2 + (k+1)l + \left\lfloor l^2/4 \right\rfloor - (k+1)^2 - (k+1)(l-1) - \left\lfloor (l-1)^2/4 \right\rfloor\notag\\
              &= k + 1 = 1 + \lfloor \omega/2 \rfloor
\end{align}

Alternatively we can easily compute $m_w$, by setting $j=\omega-i$ in the inequalities~(\ref{eq:ineq_ij_wv}) to obtain
\begin{equation}
  0 \leq \omega \leq \phi,\; 0 \leq i \leq \lfloor \omega/2 \rfloor
\end{equation}
so that from the second set of inequalities we directly obtain
\begin{equation}
  m_w(\omega) = 1 + \lfloor \omega/2 \rfloor
\end{equation}
We then calculate $n_w$, by considering separately the even and odd values of $\omega$, as
\begin{equation}
  n_w(\phi) = \sum_{\omega=0}^{\phi} \big( 1 + \lfloor \omega/2 \rfloor \big)
            = \sum_{k=0}^{\lfloor\phi/2\rfloor} ( 1 + k ) + \sum_{k=0}^{\lfloor(\phi-1)/2\rfloor} ( 1 + k )
            = \cdots = \left\lfloor \frac{(\phi+2)^2}{4} \right\rfloor
\end{equation}

The formulas for $n(\phi)$ and $m(\omega)$ for the other objective bases are obtained in a similar way. In all cases, the asymptotic behaviour for large values of $\phi$ is easily computed from the leading term of the polynomial $n(\phi)$, disregarding the presence of the floor function.

\bibliographystyle{model1-num-names}
\bibliography{Cubature}

\end{document}